\documentclass [12pt]{article}
\usepackage{latexsym}

\begin{document}

\newcommand{\Z}{\ensuremath{\mathsf{\mathbf{Z}}}}
\newcommand{\N}{\ensuremath{\mathsf{\mathbf{N}}}}
\newtheorem{thms}{Theorem}[section]
\newtheorem{lms}[thms]{Lemma}
\newtheorem{props}[thms]{Proposition}
\newtheorem{cor}[thms]{Corollary}

\begin{center}

\Large{On the independence of $K$-theory and stable rank for simple $C^{*}$-algebras}

\vspace{5mm}

\large{Andrew Toms}

\end{center}

\vspace{5mm}

\begin{abstract}

Jiang and Su and (independently) Elliott discovered a simple, nuclear, infinite-dimensional 
$C^*$-algebra $\mathcal{Z}$ having the same Elliott invariant as the
complex numbers.  For a nuclear $C^*$-algebra $A$ with weakly unperforated 
$\mathrm{K}_*$-group the Elliott invariant of $A \otimes \mathcal{Z}$ is isomorphic 
to that of $A$.  Thus, any simple nuclear $C^*$-algebra $A$ having a weakly unperforated
$\mathrm{K}_*$-group which does not absorb $\mathcal{Z}$ provides a counterexample
to Elliott's conjecture that the simple nuclear $C^*$-algebras will be classified by
the Elliott invariant.  In the sequel we exhibit
a separable, infinite-dimensional, stably finite instance of such a non-$\mathcal{Z}$-absorbing
algebra $A$, and so provide a counterexample to the Elliott conjecture
for the class of simple, nuclear, infinite-dimensional, stably finite,
separable $C^*$-algebras.

\end{abstract}

\section{Introduction}

Elliott's classification of AF $C^*$-algebras ([2]) via the scaled,
ordered $\mathrm{K}_0$-group began what is now a widespread effort to
classify nuclear $C^*$-algebras via the Elliott
invariant.  In the case of a stably finite,
unital, simple $C^*$-algebra $A$ this invariant consists of the group $\mathrm{K}_{*} A = \mathrm{K}_0 A
\oplus \mathrm{K}_1 A$, the class of the unit of $A$ in $\mathrm{K}_{*} A$,
an order structure on $\mathrm{K}_{*} A$ (an element $[p] \oplus x$ is 
positive if $[p]$ is positive in $\mathrm{K}_0 A$ and $x$ can be 
represented as the $\mathrm{K}_1$-class of a unitary $u \in \mathrm{M}_l (A)$ such that $u u^*$ is
a sub-projection of $p$), the Choquet simplex of normalised traces $\mathrm{T} A$,
and the pairing between $\mathrm{K}_0 A$ and $\mathrm{T} A$ via evaluation.  In this
paper the invariant above will be denoted $\mathrm{Ell}(A)$.  Let $\mathrm{sr}(A)$ 
be the stable rank of $A$, as defined by Rieffel in [9].
$\mathrm{Ell}(-)$ has been particularly successful in classifying simple
$C^*$-algebras of stable rank one.  Until now, it was not
known whether this invariant would suffice for the classification of
stably finite $C^*$-algebras of stable rank greater than one.

An ordered group $(G,G^+)$ is said to be weakly unperforated if 
$x \notin G^+$ and $nx \in G^+$ for some natural number $n$ implies that $nx = 0$. 
We recall that the Elliott invariant of a simple nuclear unital
$C^*$-algebra $A$ is isomorphic to that of $A \otimes \mathcal{Z}$ whenever
$\mathrm{K}_{*} A$ is weakly unperforated ([4]).  If $A \cong A \otimes 
\mathcal{Z}$, then we say that $A$ is $\mathcal{Z}$-stable.  Our main result is
the following:
\begin{thms}
For each natural number $n > 2$ there exists a simple, unital,
nuclear, separable,
infinite-dimensional, stably finite, non-$\mathcal{Z}$-stable
$C^*$-algebra $B_n$ such that $\mathrm{K}_* B_n$ is weakly
unperforated and $\mathrm{sr}(B_n) \in \{n+1,n+2\}$.  In particular,
\begin{displaymath}
\mathrm{Ell}(B_n) \simeq \mathrm{Ell}(B_n \otimes \mathcal{Z}).
\end{displaymath}
\end{thms}
Thus, $B_n$ and $B_n \otimes \mathcal{Z}$ constitute a counterexample
to the Elliott conjecture for the class of simple, nuclear,
infinite-dimensional, stably finite $C^*$-algebras.
We note that the existence of $B_n$ answers Question 1.5 of
[4] negatively;  the weak unperforation of the $\mathrm{K}_*$-group
does not imply that a simple, unital, nuclear, separable, infinite-dimensional
$C^*$-algebra absorbs $\mathcal{Z}$.

The title of this paper derives from the fact that the algebra $B_n$ of Theorem 1.1
has $\mathrm{sr}(B_n) \in \{n+1,n+2\}$ while, as we shall see, 
$\mathrm{sr}(B_n \otimes \mathcal{Z}) \leq 2$.
It is possible (but purely speculative) that finer  
invariants such as $\mathrm{K}$-theory with coefficients, the 
semigroup of Murray-von Neumann equivalence classes
of projections, or higher algebraic K-theory will recover
stable rank, and so the independence of the title is only with respect to
the notion of $\mathrm{K}$-theory captured by $\mathrm{Ell}( - )$.

We conclude this section with an outline of the sequel.
Section 2 lists several theorems from [3], which are 
applied in section 3 to construct the algebra $B_n$ of Theorem 1.1.  
The general ideas of this latter section are also found in [3].  
In section 4 $B_n$ is shown to have the properties claimed
in Theorem 1.1.

\emph{Acknowledgements.} The author would like to thank George Elliott
and Ping Wong Ng for their many helpful comments on earlier drafts of
this paper.  This research was funded by both the Israel Halperin
Graduate Award at the University of Toronto and by a Natural Sciences and
Engineering Research Council of Canada Postdoctoral Fellowship.

\section{Background and essential results}

We begin by reviewing the definition of the generalised mapping torus, due to 
Elliott.  Let $C$, $D$ be $C^{*}$-algebras and let $\phi_{0}$, $\phi_{1}$ be
$*$-homomorphisms from $C$ to $D$.  Then the generalised mapping torus of
$C$ and $D$ with respect to $\phi_{0}$ and $\phi_{1}$ is 
\begin{displaymath}
A:=\{(c,d)|d \in C([0,1];D), \ c \in C, \ d(0)=\phi_{0}(c), \ d(1)=\phi_{1}(c)\}
\end{displaymath}
We will denote $A$ by $A(C,D,\phi_{0}, \phi_{1})$ where appropriate for clarity.  
We now list (without proof) some theorems of [3] which will be used in the sequel.

\begin{thms}[{Elliott and Villadsen ([3]), Theorem 2}]
The index map $b_{*} : K_{*}C \longrightarrow K_{1-*}SD = K_{*}D$ in the six term periodic 
sequence for the extension
\begin{displaymath}
0 \longrightarrow SD \longrightarrow A \longrightarrow C \longrightarrow 0
\end{displaymath}
is the difference
\begin{displaymath}
\mathrm{K}_{*} \phi_{1} - \mathrm{K}_{*} \phi_{0} : \mathrm{K}_{*} C \longrightarrow \mathrm{K}_{*} D.
\end{displaymath}
Thus, the six-term exact sequence may be written as the short exact sequence
\begin{displaymath}
0 \longrightarrow \mathrm{Coker} b_{1-*} \longrightarrow \mathrm{K}_{*} A \longrightarrow \mathrm{Ker} 
b_{*} \longrightarrow 0.
\end{displaymath}
In particular, if $b_{1-*}$ is surjective, then $\mathrm{K}_{*} A$ is isomorphic to its image, 
$\mathrm{Ker} b_{*}$, in $\mathrm{K}_{*} C$.

Suppose that cancellation holds for $D$ --- i.e., that cancellation holds in the semigroup of 
Murray-von Neumann equivalence classes of projections of projections in $D$ and
in matrix algebras over $D$ (equivalently, in $D \otimes \mathcal{K}$).  It follows that if $b_{1}$ is 
surjective, so that $\mathrm{K}_{0} A \subseteq \mathrm{K}_{0} C$, then
\begin{displaymath}
{(\mathrm{K}_{0} A)}^{+} = {(\mathrm{K}_{0} C)}^{+} \cap \mathrm{K}_{0} A.
\end{displaymath}

The preceding conclusion also holds if cancellation is only known to hold for each pair of 
projections in $D \otimes \mathcal{K}$ obtained as the images under the maps $\phi_{0}$ and
$\phi_{1}$ of a single projection in $C \otimes \mathcal{K}$.  (In other words, if two such such 
projections in $D \otimes \mathcal{K}$ have the same $\mathrm{K}_{0}$ class then they should be 
equivalent, assuming as before that $b_{1}$ is surjective.)  
\end{thms}

\begin{thms}[{Elliott and Villadsen ([3]), Theorem 3}]

Let $A_{1}$ and $A_{2}$ be building block algebras as described above,

\begin{displaymath} A_{i}=A(C,D,\phi_{0}^{i},\phi_{1}^{i}),\ \ i=1,2.
\end{displaymath}

Let there be given four maps between the fibres,
\begin{displaymath}
\begin{array}{rrr}
\gamma : & C_{1} \rightarrow C_{2}, &  \\
\delta, \delta ' : & D_{1} \rightarrow D_{2}, & \ \ and, \\
\epsilon : & C_{1} \rightarrow D_{2}, & 
\end{array}
\end{displaymath}
such that $\delta$, $\delta '$ and $\epsilon$ have mutually orthogonal images, and
\begin{displaymath} \delta \phi_{0}^{1} + \delta ' \phi_{1}^{1} + \epsilon = \phi_{0}^{2} \gamma,
\end{displaymath}
\begin{displaymath} \delta \phi_{1}^{1} + \delta ' \phi_{0}^{1} + \epsilon = \phi_{1}^{2} \gamma.
\end{displaymath}

Then there exists a unique map 
\begin{displaymath} \theta : A_{1} \rightarrow A_{2},
\end{displaymath}
respecting the canonical ideals, giving rise to the map $\gamma : C_{1} \rightarrow
C_{2}$ between the quotients (or fibres at infinity), and such that for any $0 < s< 1$,
if $e_{s}$ denotes evaluation at $s$, and $e_{\infty}$ the evaluation at infinity,
\begin{displaymath} e_{s} \theta = \delta e_{s} + \delta ' e_{1-s} + \epsilon e_{\infty}.
\end{displaymath}

\end{thms}

\begin{thms}[{Elliott and Villadsen ([3]), Theorem 4}]

Let $A_{1}$ and $A_{2}$ be building block algebras as in Theorem 2.1.  Let $\theta : A_{1} \rightarrow 
A_{2}$ be a homomorphism constructed as in Theorem 2.2, from maps $\gamma : C_{1} \rightarrow C_{2}$, 
$\delta,\delta ' : D_{1} \rightarrow D_{2}$, and $\epsilon : C_{1} \rightarrow D_{2}$.

Let there be given a map $\beta : D_{1} \rightarrow C_{2}$ such that the composed map
$\beta \phi_{1}^{1}$ is a direct summand of the map $\gamma$, and such that the 
composed maps $\phi_{0}^{2} \beta$ and $\phi_{1}^{2} \beta$ are direct summands of
the maps $\delta '$ and $\delta$, respectively.  Suppose that the decomposition of
$\gamma$ as the orthogonal sum of $\beta \phi_{1}^{1}$ and another map is such that
the image of the second map is orthogonal to the image of $\beta$.  (Note that this
requirement is automatically satisfied if $C_{1}$, $D_{1}$, and the map $\beta \phi_{1}^{1}$
are unital.)

It follows that, for any $0 < t < \frac{1}{2}$, the map $\theta : A_{1} \rightarrow A_{2}$
is homotopic to a map $\theta_{t} : A_{1} \rightarrow A_{2}$ differing from it only as follows: 
the map $e_{\infty} \theta_{t}$ has the direct summand $\beta e_{t}$ instead of one of the 
direct summands $\beta \phi_{0}^{1} e_{\infty}$ and $\beta \phi_{1}^{1} e_{\infty}$ of 
$e_{\infty} \theta$, and for each $0 < s < 1$ the map $e_{s} \theta_{t}$ has either the
direct summand $\phi_{0}^{2} \beta e_{t}$ instead of the direct summand $\phi_{0}^{2} \beta e_{s}$
of $e_{s} \theta$, or the direct summand $\phi_{1}^{2} \beta e_{t}$ instead of the direct summand 
$\phi_{1}^{2} \beta e_{s}$ of $e_{s} \theta$, or both.

Furthermore, let $\alpha : D_{1} \rightarrow C_{2}$ be any map homotopic to $\beta$ within the 
hereditary sub-$C^{*}$-algebra of $C_{2}$ generated by the image of $\beta$.  Then the map 
$\theta_{t}$ is homotopic to a map $\theta_{t}^{'} : A_{1} \rightarrow A_{2}$ differing from $\theta_{t}$
only in the direct summands mentioned, and such that $e_{\infty} \theta_{t}^{'}$ has the direct summand
$\alpha e_{t}$ instead of $\beta e_{t}$, and for each $0 < s < 1$, $e_{s} \theta_{t}^{'}$ has
either $\phi_{0}^{2} \alpha e_{t}$ instead of $\phi_{0}^{2} \beta e_{t}$, or $\phi_{1}^{2} \alpha e_{t}$
instead of $\phi_{1}^{2} \beta e_{t}$.  

\end{thms}

\begin{thms}[{Elliott and Villadsen ([3]), Theorem 5}]

Let
\begin{displaymath} A_{1} \stackrel{\theta_{1}}{\rightarrow} A_{2} \stackrel{\theta_{2}}{\rightarrow} \cdots
\end{displaymath}
be a sequence of separable building block $C^{*}$-algebras,
\begin{displaymath} A_{i} = A(C_{i},D_{i},\phi_{0}^{i}, \phi_{1}^{i}), \ i=1,2, \ldots
\end{displaymath} 
with each map $\theta_{i} : A_{i} \rightarrow A_{i+1}$ obtained by the construction of Theorem 2.2 (and thus 
respecting the canonical ideals).  For each $i=1,2, \ldots$ let $\beta_{i} : D_{i} \rightarrow C_{i+1}$ be a
map verifying the hypotheses of Theorem 2.3.

Suppose that for every $i=1,2, \ldots$, the intersection of the kernels of the boundary maps $\phi_{0}^{i}$ and 
$\phi_{1}^{i}$ from $C_{i}$ to $D_{i}$ is zero.

Suppose that, for each i, the image of each of $\phi_{0}^{i+1}$ and $\phi_{1}^{i+1}$ generates $D_{i+1}$ as a
closed two-sided ideal, and that this is in fact true for the restriction of $\phi_{0}^{i+1}$ and 
$\phi_{1}^{i+1}$ to the smallest direct summand of $C_{i+1}$ containing the image of $\beta_{i}$.  Suppose
that the closed two-sided ideal of $C_{i+1}$ generated by the image of $\beta_{i}$ is a direct summand. 

Suppose that, for each i, the maps $\delta_{i}^{'} - \phi_{0}^{i} \beta_{i}$ and 
$\delta_{i} - \phi_{1}^{i} \beta_{i}$ from $D_{i}$ to $D_{i+1}$ are injective.  

Suppose that, for each i, the map $\gamma_{i} - \beta_{i} \phi_{1}^{i}$ takes each non-zero direct summand
of $C_{i}$ into a subalgebra of $C_{i+1}$ not contained in any proper closed two-sided ideal. 

Suppose that, for each i, the map $\beta_{i} : D_{i} \rightarrow C_{i+1}$ can be deformed---inside the 
hereditary sub-$C^{*}$-algebra generated by its image---to a map $\alpha_{i} : D_{i} \rightarrow C_{i+1}$
with the following property:  There is a direct summand of $\alpha_{i}$, say ${\bar{\alpha}}_{i}$, such
that ${\bar{\alpha}}_{i}$ is non-zero on an arbitrary given element $x_{i}$ of $D_{i}$, and has image a 
simple sub-$C^{*}$-algebra of $C_{i+1}$, the closed two-sided ideal generated by which contains the image 
of $\beta_{i}$.

Choose a dense sequence $(t_{n})$ in the open interval $(0,\frac{1}{2})$, such that $t_{2n} = t_{2n-1},
n = 1,2,\ldots$.

Choose a sequence of elements $x_{3} \in D_{3}, \ x_{5} \in D_{5}, \ x_{7} \in D_{7}, \ldots$ 
(necessarily non-zero) with the following property:  For some countable basis for the topology of 
the spectrum of each of $D_{1}, D_{2},\ldots$, and for some choice of non-zero element of the closed 
two-sided ideal associated to each of these (non-empty) open sets, under successive application of the
maps $\delta_{i} - \phi_{1}^{i+1} \beta_{i}$ each one of these elements is taken into $x_{j}$ for all j 
in some set $S \subseteq \{3,5,7,\ldots\}$ such that $\{t_{j}, j \in S \}$ is dense in $(0,\frac{1}{2})$.  
Choose $\alpha_{j}$ as above such that ${\bar{\alpha}}_{j}(x_{j}) \neq 0$ for some direct summand
${\bar{\alpha}}_{j}$ of $\alpha_{j}$ for each $j \in \{3,5,7,\ldots\}$.  For each $j \in \{4,6,8,\ldots\}$
choose $\alpha_{j}$ with respect to the non-zero element $(\delta_{j-1}^{'} - \phi_{0}^{j} \beta_{j-1})
(x_{j-1})$ of $D_{j}$.  (If $j=1$ or $2$, choose $\alpha_{j} = \beta_{j}$.)

It follows that, if $\theta_{i}^{'}$ denotes the deformation of $\theta_{i}$ constructed in Theorem 4, 
with respect to the point $t_{i} \in (0,\frac{1}{2})$ and the maps $\alpha_{i}$ and $\beta_{i}$ (and
a fixed homotopy of $\beta_{i}$ to $\alpha_{i}$), then the inductive limit of the sequence
\begin{displaymath} A_{1} \stackrel{\theta_{1}^{'}}{\rightarrow} A_{2} \stackrel{\theta_{2}^{'}}{\rightarrow}
\cdots
\end{displaymath}
is simple.

\end{thms}

\section{The construction of $B_n$}

We now specify $C^{*}$-algebras $A_{i} = A_{i}(C_{i}, D_{i}, \phi_{i}^{0}, \phi_{i}^{1})$ as in 
Theorem 2.1, and maps $\delta_{i}$, $\delta_{i}^{'}$, $\gamma_{i}$, and $\beta_{i}$ satisfying 
the hypotheses of Theorems 2.2, 2.3, and 2.4 in order to construct an inductive sequence
\begin{displaymath}
A_{1} \stackrel{\theta_{1}^{'}}{\longrightarrow} A_{2} \stackrel{\theta_{1}^{'}}{\longrightarrow} \cdots
\end{displaymath}
whose limit will be the $C^*$-algebra $B_n$ of Theorem 1.1. 

Let $\mathrm{D}$ denote the closed unit disc in the complex numbers. Put
\begin{displaymath}
X_{i} = \mathrm{D}^{n} \times \mathrm{CP}^{n\sigma(1)} \times
\mathrm{CP}^{n\sigma(2)} \times \cdots \times \mathrm{CP}^{n\sigma(i)}
\end{displaymath}
--- the $\sigma(i)$ are natural numbers to be specified --- so that 
\begin{displaymath}
X_{i+1} = X_{i} \times \mathrm{CP}^{n\sigma(i+1)}, 
\end{displaymath}
and let
\begin{displaymath}
\begin{array}{cc}
\pi_{i+1}^{1}: X_{i+1} \longrightarrow X_{i}, \  & \ \pi_{i+1}^{2}:
X_{i+1} \longrightarrow \mathrm{CP}^{n\sigma(i+1)}
\end{array}
\end{displaymath}
be the co-ordinate projections.

We will take $C_{i} = p_{i}(\mathrm{C}(X_{i}) \otimes \mathcal{K})p_{i}$, where $p_{i}$ is a projection in
$\mathrm{C}(X_{i}) \otimes \mathcal{K}$ to be specified.  Let $D_{i} =
C_{i} \otimes  \mathrm{M}_{k_{i} \mathrm{dim}(p_{i})}$,
where $k_{i}$ is a positive integer to be specified.  Define maps 
\begin{displaymath}
\mu_{i}, \nu_{i}:C_{i} \longrightarrow C_{i} \otimes \mathrm{M}_{\mathrm{dim}(p_{i})}
\end{displaymath}
as follows:
\begin{displaymath}
\mu_{i}(a) = p_{i} \otimes a(x_{i}),
\end{displaymath}
\begin{displaymath}
\nu_{i}(a) = a \otimes 1_{\mathrm{dim}(p_{i})}.
\end{displaymath}
For $t \in \{0,1\}$, we will take $\phi_{i}^{t}$ to be the direct sum of $l_{i}^{t}$ copies of 
$\mu_{i}$ and $k_{i} - l_{i}^{t}$ copies of $\nu_{i}$, where the $l_{i}^{t}$ are non-negative integers
to be specified.  All that we mention now is that we should have $l_{i}^{1} - l_{i}^{0} \neq 0$.
We need only specify the $\phi_{i}^{t}$ up to unitary equivalence, a fact we shall
exploit below.

By Theorem 2.1 we have that for any $e \in \mathrm{K}_{0}(C_{i})$,
\begin{displaymath}
\begin{array}{rl}
b_{0}(e) & = (l_{i}^{1} - l_{i}^{0})(\mathrm{K}_{0}(\mu_{i}) - \mathrm{K}_{0}(\nu_{i})) \\
         & = (l_{i}^{1} - l_{i}^{0})(\mathrm{dim}(e) \cdot \mathrm{K}_{0}(p_{i}) - \mathrm{dim}(p_{i}) \cdot e).
\end{array}
\end{displaymath}
Since $l_{i}^{1} - l_{i}^{0} \neq 0$ and since $\mathrm{K}_{0}C_{i}$ is a finitely generated free abelian group, 
we have that $\mathrm{Ker}b_{0}$ is the largest subgroup 
of $\mathrm{K}_{0}C_{i}$ containing $\mathrm{K}_{0}(p_{i})$ and
isomorphic to the integers.  In the sequel we will choose $p_{i}$
so that $\mathrm{K}_{0}(p_{i})$ in fact generates said subgroup.  
Since $\mathrm{K}_{1}C_{i} = 0$ we have, by Theorem
2.1, that $\mathrm{K}_{0}A_{i}$ is isomorphic as an ordered group to
its image, $\mathrm{Ker}b_{0}$, in $\mathrm{K}_{0}C_{i}$,
considered as a sub ordered group.  The latter (with the choice of $p_i$
below) is isomorphic to the integers with the unique unperforated
order structure, and the image of $[1_{A_i}]$ is $[p_i]$.

Let $p_{1}$ be a projection corresponding to the vector bundle 
\begin{displaymath}
\theta_{1} \times \xi_{n\sigma(1)},
\end{displaymath}
over $X_{1}$, where $\theta_{1}$ denotes the trivial line bundle of
dimension one over $\mathrm{D}$, $\xi_k$
denotes the universal line bundle over $\mathrm{CP}^{k}$ for every
natural number $k$, and $\sigma(1) = 1$.  We now 
specify, inductively, the maps $\gamma_{i}:C_{i} \longrightarrow C_{i+1}$.  Consider first the map
\begin{displaymath}
\psi_{i} := \mathrm{id} \otimes 1
\end{displaymath}
from $\mathrm{C}(X_{i})$ to $\mathrm{C}(X_{i+1}) = \mathrm{C}(X_{i}
\times \mathrm{CP}^{n\sigma(i+1)}) = \mathrm{C}(X_{i}) \otimes \mathrm{C}(\mathrm{{}CP}^{n\sigma(i+1)})$, 
where $1$ denotes the unit of $\mathrm{C}(\mathrm{CP}^{n\sigma(i+1)})$ and $\mathrm{id}$ denotes the identity map from
$\mathrm{C}(X_{i})$ to itself.

Consider also the map
\begin{displaymath}
\beta_{i}^{'} := {\pi_{i+1}^{2 *}}(\xi_{n\sigma(i+1)}) \cdot e_{x_{i}}
\end{displaymath}
from $\mathrm{C}(X_{i})$ to $\mathrm{C}(X_{i+1}) \otimes \mathcal{K}$, where 
$e_{x_{i}}$ denotes evaluation at $x_{i}$.  All that we shall require of the
$x_{i}$ at this stage is that $\pi_{i+1}^{1}(x_{i+1}) = x_{i}$.

Now, inductively, let us take $\gamma_{i}$ to be the map from $C_{i}$ 
to $\mathrm{C}(X_{i+1}) \otimes \mathrm{M}_{2}(\mathcal{K})$ consisting
of the direct sum of the following two maps:  first, the restriction 
to $C_{i} \subseteq \mathrm{C}(X_{i}) \otimes \mathcal{K}$ of the
tensor product of $\psi_{i}$ with the identity map from $\mathcal{K}$ 
to $\mathcal{K}$, and second, the map from $C_{i}$ to
$\mathrm{C}(X_{i+1}) \otimes \mathrm{M}_{q_{i}}(\mathcal{K})$ 
consisting of the composition of the map $\phi_{i}^{1}$ from $C_{i}$ to 
$D_{i}$ with the direct sum of $q_{i}$ copies of the tensor product of the map $\beta_{i}^{'}$ with the 
identity map from $\mathcal{K}$ to $\mathcal{K}$ (restricted to $D_{i}
\subseteq \mathrm{C}(X_{i}) \otimes \mathcal{K}$), where $q_{i}$ is to
be specified.  The induction consists in first considering the case 
$i=1$ (as $p_{1}$ has already been chosen), 
then setting $p_{2} = \gamma_{1}(p_{1})$, so that $C_{2}$ is 
specified as the cut-down of $\mathrm{C}(X_{2}) \otimes \mathcal{K}$
by $p_{2}$, and continuing in this way.

With the maps $\gamma_{i}$ defined as above, we have that $p_{i}$ 
is a projection in $\mathrm{C}(X_{i}) \otimes \mathcal{K}$
corresponding to the vector bundle
\begin{displaymath}
\theta_{1} \times \xi_{n} \times \sigma(2) \xi_{n\sigma(2)} \times \cdots \times \sigma(i) \xi_{n\sigma(i)},
\end{displaymath}
where 
\begin{displaymath}
\sigma(i) = \prod_{l=1}^{i-1} (\mathrm{mult}(\gamma_{l}) - 1).
\end{displaymath}
Notice that by the K{\"{u}}nneth formula (in Chapter 5 of [10], for instance) the classes 
$[\theta_1]$, $[\xi_{n \sigma(1)}]$, \ldots, $[\xi_{n \sigma(i)}]$ are independent in
$\mathrm{K}^{0}(X_i)$ (we are abusing notation slightly here, using $[\xi_k]$
to represent the class of the induced bundle $\pi^{*}(\xi_k)$, where $\pi$ is  
projection from $X_i$ onto $\mathrm{CP}^k$).  Suppose that $[p_{i}] = ky$ for some $k \in
\mathrm{Z}$, $y \in \mathrm{K}^0(X_i)$.  It follows from independence
that we have $[\xi_n] = ky^{'}$, $y^{'} \in
\mathrm{K}^0(\mathrm{CP}^n)$, whence $k = \pm 1$.  We conclude that
$[p_i]$ itself generates the subgroup of rational multiples of
$[p_i]$ in $\mathrm{K}^{0}X_{i}$, as desired.  Thus $\gamma_{i}$ 
induces an isomorphism of ordered groups from
$\mathrm{Ker}b_{0}$ at the $i^{\mathrm{th}}$ stage 
to $\mathrm{Ker}b_{0}$ at the ${(i+1)}^{\mathrm{th}}$ stage.

Note that $\gamma_{i} - \beta_{i} \phi_{i}^{1}$ is non-zero, and so takes $C_{i}$ into a subalgebra of $C_{i+1}$
not contained in any proper closed two-sided ideal.

Next, we construct the maps $\delta_{i}, \delta_{i}^{'}:D_{i} \rightarrow D_{i+1}$, with orthogonal images,
such that
\begin{displaymath}
\begin{array}{c}
\delta_{i} \phi_{i}^{0} + \delta_{i}^{'} \phi_{i}^{1} = \phi_{i+1}^{0} \gamma_{i}, \\
\delta_{i} \phi_{i}^{1} + \delta_{i}^{'} \phi_{i}^{0} = \phi_{i+1}^{1} \gamma_{i},
\end{array}
\end{displaymath}
and $\phi_{i+1}^{0} \beta_{i}$ and $\phi_{i+1}^{1} \beta_{i}$ are direct summands of $\delta_{i}^{'}$ and 
$\delta_{i}$, respectively.  To do this we shall have to modify $\phi_{i+1}^{0}$ and $\phi_{i+1}^{1}$ by 
inner automorphisms.  This will not affect the $\mathrm{K}$-theory of $A$.

Note that, up to unitary equivalence, we have
\begin{displaymath}
e_{x_{i+1}} \gamma_{i} = \mathrm{mult}(\gamma_{i}) e_{x_{i}},
\end{displaymath}
where $\mathrm{mult}(\gamma_{i})$ denotes the factor by which
$\gamma_{i}$ multiplies dimension.  It follows that up to unitary equivalence
\begin{displaymath}
\begin{array}{rl}
\mu_{i+1} \gamma_{i} &  = p_{i+1} \otimes e_{x_{i+1}} \gamma_{i} \\
                     &  = \gamma_{i}(p_{i}) \otimes \mathrm{mult}(\gamma_{i}) e_{x_{i}} \\
                     &  = \mathrm{mult}(\gamma_{i}) \gamma_{i}(p_{i} \otimes e_{x_{i}}) \\
                     &  = \mathrm{mult}(\gamma_{i}) \gamma_{i} \mu_{i},
\end{array}
\end{displaymath}
and
\begin{displaymath}
\begin{array}{rl}
\nu_{i+1} \gamma_{i} &  = \gamma_{i} \otimes 1_{\mathrm{dim}(p_{i+1})} \\
                     &  = \mathrm{mult}(\gamma_{i}) \gamma_{i} \otimes 1_{\mathrm{dim}(p_{i})} \\
                     &  = \mathrm{mult}(\gamma_{i}) \gamma_{i} \nu_{i}.
\end{array}
\end{displaymath}

Take $\delta_{i}$ and $\delta_{i}^{'}$ to be $r_{i}$ and $s_{i}$ copies of $\gamma_{i}$, where $r_{i}$ and
$s_{i}$ are integers to be specified.  The conditions
\begin{displaymath}
\delta_{i} \phi_{i}^{0} + \delta_{i}^{'} \phi_{i}^{1} = \phi_{i+1}^{0} \gamma_{i}
\end{displaymath}
and
\begin{displaymath}
\delta_{i} \phi_{i}^{1} + \delta_{i}^{'} \phi_{i}^{0} = \phi_{i+1}^{1} \gamma_{i},
\end{displaymath}
understood up to unitary equivalence imply that
\begin{displaymath}
r_{i} \gamma_{i}(l_{i}^{t}\mu_{i} - (k_{i} - l_{i}^{t})\nu_{i}) + s_{i} \gamma_{i}(l_{i}^{1-t}\mu_{i} + 
(k_{i} - l_{i}^{1-t})\nu_{i}) = (l_{i+1}^{t}\mu_{i+1} + (k_{i+1} - l_{i+1}^{t})\nu_{i+1}) \gamma_{i},
\end{displaymath}
again, up to unitary equivalence.  As $\mathrm{K}_{0}(\mu_{i})$ 
and $\mathrm{K}_{0}(\nu_{i})$ are independent, the above equation is
equivalent to the two equations
\begin{displaymath}
\begin{array}{rl}
r_{i} l_{i}^{t} + s_{i} l_{i}^{1-t} & = \mathrm{mult}(\gamma_{i}) l_{i+1}^{t}, \\
(r_{i} + s_{i}) k_{i} & = \mathrm{mult}(\gamma_{i}) k_{i+1}.
\end{array}
\end{displaymath}
Choose $r_{i} = 2 \mathrm{mult}(\gamma_{i})$ and $s_{i} = \mathrm{mult}(\gamma_{i})$, so that
\begin{displaymath}
k_{i+1} = 3 k_{i},
\end{displaymath}
and 
\begin{displaymath}
l_{i+1}^{t} = 2 l_{i}^{t} + l_{i}^{1-t}.
\end{displaymath}

Take $k_{1} = 1$, $l_{1}^{1} = 1$, and $l_{1}^{0} = 0$.  Then $l_{i}^{1} - l_{i}^{0} \neq 0$ for
all $i$, as required.

Next, let us show that, up to unitary equivalence preserving the equations $\delta_{i} \phi_{i}^{t} +
\delta_{i}^{'} \phi_{i}^{1-t} = \phi_{i+1}^{t} \gamma_{i}$, $\phi_{i+1}^{1} \beta_{i}$ is a direct summand
of $\delta_{i} = 2 \mathrm{mult}(\gamma_{i}) \gamma_{i}$ and $\phi_{i+1}^{0} \beta_{i}$ is a direct summand of
$\delta_{i}^{'} = \mathrm{mult}(\gamma_{i}) \gamma_{i}$.

Note that $\phi_{i+1}^{t} \beta_{i}$ is a direct sum of $l_{i+1}^{t}$ copies of $p_{i+1} \otimes \beta_{i}$
and \mbox{${(k_{i+1} - l_{i+1}^{t})}$} copies of $\beta_{i}$, whereas $\delta_{i}$ and $\delta_{i}^{'}$ contain,
respectively, $q_{i} \mathrm{mult}(\gamma_{i})$ and $2 q_{i} \mathrm{mult}(\gamma_{i})$ copies of $\beta_{i}$.
Note also, that by Theorem 8.1.2 of [5], a trivial projection of dimension $\mathrm{dim}(p_{i+1}) + \frac{1}{2}
\mathrm{dim}X_{i+1}$ in $\mathrm{C}(X_{i+1}) \otimes \mathcal{K}$ contains a copy of
$p_{i+1}$.  Therefore, $2 \mathrm{dim}(p_{i+1}) + 2
\mathrm{dim}X_{i+1}$ copies of $\beta_{i}$ contain a copy of $p_{i+1} \otimes \beta_{i}$ 
(since $2 \mathrm{dim}(p_{i+1}) + 2\mathrm{dim}X_{i+1}$ copies of $\xi_{i+1}$ contain
a trivial projection of dimension $\mathrm{dim}(p_{i+1}) + \frac{1}{2}
\mathrm{dim}X_{i+1}$).  It follows that
$k_{i+1} (2 \mathrm{dim}(p_{i+1}) + 2 \mathrm{dim} X_{i+1})$ copies of $\beta_{i}$ contain a copy of 
$\phi_{i+1}^{t} \beta_{i}$ when $t$ is equal to either $0$ or $1$.  By a copy of a given map from $D_{i}$ to
$D_{i+1}$ we mean another map obtained from it by conjugating by a partial isometry in $D_{i+1}$ with initial
projection the image of the unit.

Note that
\begin{displaymath}
k_{i+1}(2 \mathrm{dim} (p_{i+1}) + 2 \mathrm{dim} X_{i+1}) = 6 k_{i} \mathrm{mult}(\gamma_{i})(\mathrm{dim}(p_{i}) + 
                                                                \mathrm{dim} X_{i}), 
\end{displaymath}                                                         
and that $k_{i}$, $\mathrm{dim}(p_{i})$, and $\mathrm{dim}X_{i}$ have already been specified, and do not depend on $q_{i}$.
It follows that, with
\begin{displaymath}
q_{i} \geq 6 k_{i} (\mathrm{dim}(p_{i}) + \mathrm{dim} X_{i}),
\end{displaymath}
$q_i \mathrm{mult}(\gamma_{i})$ copies of $\beta_{i}$ contain a copy of $\phi_{i+1}^{t} \beta_{i}$, ($t \in \{0,1\}$).
In particular $\delta_{i}^{'}$ and $\delta_{i}$ contain copies of 
$\phi_{i+1}^{0} \beta_{i}$ and $\phi_{i+1}^{1} \beta_{i}$, respectively.

With $q_{i}$ as above, let us show that for each $t=0,1$ there exists a unitary $u_{t} \in D_{i+1}$ such that
\begin{displaymath}
(\mathrm{Ad} u_{t}) \phi_{i+1}^{t} \gamma_{i} = \phi_{i+1}^{t} \gamma_{i},
\end{displaymath}
with $(\mathrm{Ad} u_{0}) \phi_{i+1}^{0} \beta_{i}$ a direct summand of $\delta_{i}^{'}$ and $(\mathrm{Ad} u_{1})
\phi_{i+1}^{1} \beta_{i}$ a direct summand of $\delta_{i}$.  In 
other words, for each $t=0,1$, we must show that the 
partial isometry constructed in the preceding paragraph, producing 
a copy of $\phi_{i+1}^{t} \beta_{i}$ inside $\delta_{i}$
or $\delta_{i}^{'}$ may be chosen in such a way that it extends 
to a unitary element of $D_{i+1}$ --- which in addition
commutes with the image of $\phi_{i+1}^{t} \gamma_{i}$.

        Let us consider the case $t=0$; the case $t=1$ is 
similar.  First note that the unit of the image of 
$\phi_{i+1}^{0} \beta_{i}$ --- the initial projection of 
the partial isometry transforming $\phi_{i+1}^{0} \beta_{i}$
into a direct summand of $\delta_{i}^{'}$ --- lies in the 
commutant of the image of $\phi_{i+1}^{0} \gamma_{i}$.
Indeed, this projection is the image by $\phi_{i+1}^{0} \beta_{i}$ 
of the unit of $D_{i}$, which, by construction,
is the image of the unit of $C_{i}$ by $\phi_{i}^{1}$.  The 
property that $\beta_{i} \phi_{i}^{1}$ is a direct summand
of $\gamma_{i}$ implies that the image by $\beta_{i} \phi_{i}^{1}$ 
of the unit of $C_{i}$ commutes with the image of
$\gamma_{i}$.  The unit of the image of $\phi_{i+1}^{0} \beta_{i}$ 
therefore commutes with the image of 
$\phi_{i+1}^{0} \gamma_{i}$, as desired.

        The final projection of the above partial isometry also 
commutes with the image of $\phi_{i+1}^{0} \gamma_{i}$.
Indeed, it is the unit of the image of a direct summand of 
$\delta_{i}^{'}$, and since $D_{i}$ is unital it is the image
of the unit of $D_{i}$ by this direct summand; since $C_{i}$ is 
unital and $\phi_{i}^{1}: C_{i} \rightarrow D_{i}$ is 
unital, the projection in question is the image of the unit of $C_{i}$
by a direct summand of $\delta_{i}^{'} \phi_{i}^{1}$.
But $\delta_{i}^{'} \phi_{i}^{1}$ is itself a direct summand 
of $\phi_{i+1}^{0} \gamma_{i}$, and so the projection in
question is the image of the unit of $C_{i}$ by a direct 
summand of $\phi_{i+1}^{0} \gamma_{i}$, and in particular
commutes with the image of $\phi_{i+1}^{0} \gamma_{i}$.

Note that both direct summands of $\phi_{i+1}^{0} \gamma_{i}$ 
under consideration ($\phi_{i+1}^{0} \beta_{i}
\phi_{i}^{1}$ and a copy of it) factor through the evaluation 
of $C_{i}$ at the point $x_{i}$, and so are contained in the
largest such direct summand of $\phi_{i+1}^{0} \gamma_{i}$;  
this largest direct summand, say $\pi_{i}$, is seen to exist
by inspection of the construction of $\phi_{i+1}^{0} \gamma_{i}$.  
Since both projections under consideration (the images
of $1 \in C_{i}$ by the two copies of $\phi_{i+1}^{0} \beta_{i} 
\phi_{i}^{1}$) are less than $\pi_{i}(1)$, to show that they
are unitarily equivalent in the commutant of the image of 
$\phi_{i+1}^{0} \gamma_{i}$ (in $D_{i+1}$) it is sufficient to 
show that they are unitarily equivalent in the commutant of the 
image of $\pi_{i}$ in $\pi_{i}(1) D_{i+1} \pi_{i}(1)$.
Note that this image is isomorphic to $\mathrm{M}_{\mathrm{dim}
(p_{i})}(\mathrm{C})$.  By construction, the two projections in question are
Murray-von Neumann equivalent in $D_{i+1}$ and hence in $\pi_{i}(1) 
D_{i+1} \pi_{i}(1)$, but all we shall use from this is
that they have the same class in $\mathrm{K}^{0} X_{i+1}$.  
Note that the dimension of these projections is $(k_{i+1} 
\mathrm{dim}(p_{i+1}))(k_{i} \mathrm{dim}(p_{i}))$, and that 
the dimension of $\pi_{i}(1)$ is at least $l_{i+1}^{0}
{(\mathrm{dim}(p_{i+1}))}^{2}$.  Since the two projections 
under consideration commute with $\pi_{i}(C_{i})$, and this is
isomorphic to $\mathrm{M}_{\mathrm{dim}(p_{i})}(\mathrm{C})$, 
to prove unitary equivalence in the commutant of $\pi_{i}(C_{i})$ in 
$\pi_{i}(1) D_{i+1} \pi_{i}(1)$ it is sufficient to prove 
unitary equivalence of the product of these projections with
a fixed minimal projection of $\pi_{i}(C_{i})$, say $e$.  Since 
$\mathrm{K}^{0} X_{i+1}$ is torsion free, the products of the two 
projections under consideration with $e$ still have the same class 
in $\mathrm{K}^{0} X_{i+1}$.  To prove that they are unitarily
equivalent in $e D_{i+1} e$, it is sufficient (and necessary) 
to prove that both they and their complements inside $e$ are
Murray-von Neumann equivalent.  Since both the cut-down projections 
and their complements inside $e$ have the same class in
$\mathrm{K}^{0} X_{i+1}$, to prove that they (i.e., the two pairs) 
are equivalent it is sufficient, by Theorem 8.1.5 of [5], to show that
all four projections have dimension at least $\frac{1}{2} 
\mathrm{dim}X_{i+1}$ (note that $\mathrm{dim} X_{i}$ is even).
Dividing the numbers above by $\mathrm{dim}(p_{i})$ (the order 
of the matrix algebra), we see that the dimension of the 
first pair of projections is $k_{i+1} k_{i} \mathrm{mult}(\gamma_{i}) 
\mathrm{dim}(p_{i})$, so that the dimension of the 
second pair of projections is at least $\mathrm{mult}(\gamma_{i})
(l_{i+1}^{0} \mathrm{dim}(p_{i+1}) - k_{i+1} k_{i} 
\mathrm{dim}(p_{i}))$.  By construction, $\mathrm{dim}(p_{i}) = 
\frac{1}{2} \mathrm{dim} X_{i}$.  Since $k_{i+1} k_{i}$ is 
non-zero for all $i$, the first inequality holds.  Since
$l_{i+1}^{0}$, the second inequality holds if $\mathrm{mult}(\gamma_{i})$
is strictly greater than $k_{i+1} k_{i}$.  Since $k_{i+1} k_{i} = 3 
{k_{i}}^{2}$, and $k_{i}$ was chosen before $q_{i}$, we
may modify our choice of $q_{i}$ to ensure that $\mathrm{mult}(\gamma_{i})$ is sufficiently large.

This shows that the two projections in $D_{i+1}$ under consideration 
are unitarily equivalent by a unitary in the commutant of 
the image of $\phi_{i+1}^{0} \gamma_{i}$.  Replacing $\phi_{i+1}^{0}$ 
by its composition with the corresponding inner
automorphism, we may suppose that the two projections in question 
are equal.  In other words $\phi_{i+1}^{0} \beta_{i}$ is
unitarily equivalent to the cut-down of $\delta_{i}^{'}$ by 
the projection $\phi_{i+1}^{0} \beta_{i}(1)$.

Now consider the compositions of these two maps with $\phi_{i}^{1}$,
i.e., $\phi_{i+1}^{0} \beta_{i} \phi_{i}^{1}$ 
and the cut-down of $\delta_{i}^{'} \phi_{i}^{1}$ by the projection
$\phi_{i+1}^{0} \beta_{i}(1)$. Since both of these maps can be 
viewed as the cut-down of $\phi_{i+1}^{0} \gamma_{i}$ by the 
same projection, they are in fact the same map.  

Therefore, any unitary inside the cut-down of $D_{i+1}$ by 
$\phi_{i+1}^{0} \beta_{i}(1)$ taking $\phi_{i+1}^{0} \beta_{i}$
into the cut-down of $\delta_{i}^{'}$ by this projection --- 
such a unitary is known to exist --- must commute with the
image of $\phi_{i+1}^{0} \beta_{i} \phi_{i}^{1}$, and hence 
with the image of $\phi_{i+1}^{0} \gamma_{i}$ --- since this 
commutes with the projection $\phi_{i+1}^{0} \beta_{i}(1)$.  
The extension of such a partial unitary to a unitary $u_{0}$ 
in $D_{i+1}$ equal to one inside the complement of this projection 
then belongs to the commutant of the image of 
$\phi_{i+1}^{0} \gamma_{i}$, and transforms $\phi_{i+1}^{0} 
\beta_{i}$ into the cut-down of $\delta_{i}^{'}$ by this
projection, as desired.

Inspection of the construction will show that the maps $\delta_{i}^{'}
- \phi_{i}^{0} \beta_{i}$ and $\delta_{i} - \phi_{i}^{1}
\beta_{i}$ are injective, as required in the hypotheses of Theorem 2.4.

Replacing $\phi_{i+1}^{t}$ with $(\mathrm{Ad} u_{t}) \phi_{i+1}^{t}$
and deforming the $\beta_i$ to other point evaluations $\alpha_i$
which are non-zero on a given element (as we may, since $X_i$ is
connected), we have completed the construction of the desired
inductive system $(A_{i}, \theta_{i})$ satisfying
the hypotheses of Theorem 2.4.  Thus, the limit $B_n$ of the 
inductive system with deformed finite stage maps, $(A_{i}, \theta_{i}^{'})$, 
is simple.  Notice
that $(\mathrm{K}_{0}B_n, [1_{B_n}]) = (\Z,1)$ --- the
$\theta_i^{'}$ are unital and  $(\mathrm{K}_{0}A_i, [1_{A_i}]) 
= (\Z,1)$ for every $i$ --- and that $B_n$ is separable, nuclear and
stably finite since each of the $A_i$ is ([1]).

\section{The main result}

In this section we prove Theorem 1.1 through a series of lemmas.  We
establish that $\mathrm{sr}(B_n) \in \{n+1,n+2\}$ (Lemma 4.1), that
$\mathrm{K}_*$ is weakly unperforated (Lemma 4.3), and that
$B_n$ does not absorb $\mathcal{Z}$ (Lemma 4.4).  Taken together,
these results show that $B_n$ is as claimed in Theorem 1.1.

\begin{lms}
\begin{displaymath}
\mathrm{sr}(B_n) \in \{n+1,n+2\}.
\end{displaymath}
\end{lms}

The proof will depend on some definitions and results which we review below.

For a unital $C^*$-algebra $A$ we let
\begin{displaymath}
\mathrm{Lg}_s(A) = \{(a_1,\ldots,a_s) \in A^s| a_1 A + \ldots + a_s A = A \}
\end{displaymath}
for every natural number $s$, and recall that the stable rank of $A$,
$\mathrm{sr}(A)$, is the 
least natural number $s$ such that $\mathrm{Lg}_s(A)$ is dense in $A^s$.  If no 
such natural number exists, we set $\mathrm{sr}(A) = \infty$ ([9]).
Note that if $(c_k,d_k)$ are elements of a generalised
mapping torus $A(C,D,\phi_0,\phi_1)$ for $k \in \{1,2,\ldots,n\}$ such
that 
\begin{displaymath}
\mathrm{dist} ((c_1,c_2,\ldots,c_n),\mathrm{Lg}_n (C)) \geq \delta, 
\end{displaymath}
then 
\begin{displaymath}
\mathrm{dist} (((c_1,d_1),(c_2,d_2),\ldots,(c_n,d_n)),\mathrm{Lg}_n
(A)) \geq \delta.  
\end{displaymath}
Indeed, one can check that 
\begin{displaymath}
\|{(c,d)}\| := \mathrm{max} \{ \|{c}\|, \mathrm{sup}_{t \in [0,1]} \|{d(t)}\| \}
\end{displaymath}
defines the unique $C^*$-norm on $A(C,D,\phi_0,\phi_1)$.  Thus, if 
$\mathrm{dist}(c,c^{'}) \geq \delta$ for $c,c^{'} \in C$, then 
$\mathrm{dist}((c,d),(c^{'},d^{'})) \geq \delta$ for any $(c,d),(c^{'},d^{'}) 
\in A(C,D,\phi_0,\phi_1)$.

For the remainder of this proof, any notation with subscript $i$ 
refers, where applicable, to the corresponding object in section 3. 
In order to show that $B_n$ has stable rank greater than $n$, we must
exhibit $n$ sequences of elements $A_i \ni a_{i,j} = \theta_{i1}^{'}(a_{1,j})
= (c_{i,j},d_{i,j})$, $1 \leq j \leq n$, $i \in \N$, such that 
\begin{displaymath}
\mathrm{dist} (((c_{i,1},d_{i,1}),(c_{i,2},d_{i,2}),\ldots,(c_{i,n},d_{i,n})), \mathrm{Lg_{n}}
(A_i)) \geq \delta > 0 
\end{displaymath}
for all $i$.  From this it follows that 
\begin{displaymath}
\mathrm{dist} ((\theta_{\infty 1}^{'}((c_{1,1},d_{1,1})),\ldots,
\theta_{\infty 1}^{'}((c_{1,n},d_{1,n}))),\mathrm{Lg}_n(B_n)) \geq
\delta,
\end{displaymath}
so that $\mathrm{sr}(B_n) > n$ by definition.  (Here $\theta_{\infty
  1}^{'}$ denotes the inclusion of $A_i$ into $B_n$.)
By the definition of the norm on the $A_i$, it will be enough to show that 
\begin{displaymath}
\mathrm{dist}
((c_{i,1},c_{i,2},\ldots,c_{i,n}), \mathrm{Lg_{n}} (C_i)) 
\geq \delta > 0
\end{displaymath}
for all $i$.  

We now review Theorem 7 of [11].  Let $e(\cdot)$ denote the Euler class
of a vector bundle.  Suppose that $C$ is a 
$C^*$-algebra of the form 
\begin{displaymath}
(r+q)(\mathrm{C}(M \times \mathrm{D}^n) 
\otimes \mathcal{K})(r+q),
\end{displaymath}
where $M$ is a smooth oriented manifold, and $r$ 
and $q$ are orthogonal projections in $\mathrm{C}(M \times \mathrm{D}^n) 
\otimes \mathcal{K}$ such that $r$ corresponds to the trivial line bundle 
and $q$ corresponds to a vector bundle $\alpha$ for which ${e(\alpha)}^n \neq 0$.
Let $\pi: M \times \mathrm{D}^n \rightarrow \mathrm{D}^n$ be
projection onto $\mathrm{D}^n$, and let $f_j: \mathrm{D}^n \rightarrow \mathrm{D}$ be the
$j^{\mathrm{th}}$ co-ordinate projection.

\begin{thms}[{Villadsen ([11]), Theorem 7}]
Let $C$, $\pi$ and $f_j$ be as above, and let $\tilde{c} = (c_1,\ldots,c_n) \in C^n$ be 
such that $r c_j r = (f_j \circ \pi)r$ for all $1 \leq j \leq n$.  Then, 
$\mathrm{dist}(\tilde{c}, \mathrm{Lg}_n (C)) \geq 1$.
\end{thms}

\noindent
\textbf{Proof of Lemma 4.1:}
We wish to apply Theorem 4.2 above to the algebras $C_i$, 
$i \geq 1$.  The sequel is similar to the proof of Theorem 8 
in [11].  For all $i$, let $r_i$ denote
the sub-projection of the unit of $C_i$ corresponding to the
one-dimensional trivial sub-bundle of $\theta_1 \times \xi_n \times
\cdots \times \sigma(i) \xi_{n \sigma(i)}$.  Note that $p_i$ considered
as a vector bundle over $X_i$ is the Whitney sum of $r_i$ and a second 
vector bundle, say $q_i$, and this second vector bundle has 
${e(q_i)}^n \neq 0$.  Indeed, 
\begin{displaymath}
q_i = \xi_{n} \times \sigma(2) \xi_{n\sigma(2)} \times \cdots \times \sigma(i) \xi_{n\sigma(i)},
\end{displaymath} 
and $e(\omega \oplus \gamma) = e(\omega) e(\gamma)$ for any two vector bundles
$\omega$ and $\gamma$ over a fixed base space so that 
\begin{displaymath}
{e(q_i)}^n = {e(\xi_n)}^n {e(\xi_{n \sigma(2)})}^{n \sigma(2)} \cdots  
{e(\xi_{n \sigma(i)})}^{n \sigma(i)}.
\end{displaymath}
(We are, as before, abusing notation slightly, using $\xi_k$
to represent the bundle induced on $X_i$ by $\xi_k$ via projection
from $X_i$ onto $\mathrm{CP}^k$.)
Since the integral cohomology ring $\mathrm{H}^* (\mathrm{CP}^k)$ is generated by
$e(\xi_k)$ with the relation ${e(\xi_k)}^{k+1}=0$, we may conclude by
the K{\"{u}}nneth Theorem that ${e(q_i)}^n \neq 0$, as claimed.  Each $X_i$ is 
of the form $M_i \times \mathrm{D}^n$ for some smooth oriented manifold
$M_i$, so the $C_i$ have the same form as the algebra $C$ of 
Theorem 4.2.

Note that for any element $c \in C_i$ there exists an element $(c,d)
\in A_i$ for some suitable $d \in \mathrm{C}([0,1];D)$.
Let $\pi_{i}:X_{i} \rightarrow D^n$ be the co-ordinate projection,
and let $f_j:D^n \rightarrow D$ be projection onto the $j^{th}$
co-ordinate.  Let $a_{1,j} = (c_{1,j},d_{1,j})$ be elements of $A_1$
such that $c_{1,j} = (f_j \circ \pi_1) r_1$, $1 \leq j \leq n$.  For each $i \geq 2$, put
$a_{i,j} = \theta_{i-1}^{'} \circ \theta_{i-2}^{'} \circ \cdots \circ
\theta_{1}^{'} (a_{1,j})$.  Write $a_{i,j} = (c_{i,j},d_{i,j})$.

In section 3, the map $\gamma_i$ was constructed as the direct sum of 
$\psi_i$ and a second map.  Let 
$\psi_{i1}$ denote the composition $\psi_i \circ \psi_{i-1} \circ 
\cdots \circ \psi_1$.  Note that $\psi_{i1}(r_1) = r_{i+1}$.  
By Theorem 2.3 
\begin{displaymath}
c_{i+1,j} = \psi_i(c_{i,j}) \oplus c_{i+1,j}^{'}, 
\end{displaymath}
where $c_{i+1,j}^{'}$ is an element
of the cut down of $C_i$ by $q_i$;  the deformation
of $\theta_i$ to $\theta_i^{'}$ is visible in the fibre at infinity
only in the perturbation of the image of the second direct summand of
$\gamma_i$ --- the image of $\psi_i$ remains unchanged.  Thus, by construction 
\begin{displaymath}
r_{i+1} c_{i+1,j} r_{i+1} = \psi_{i1}(r_1) c_{i+1,j} \psi_{i1}(r_1) =
\psi_{i1}(c_{1,j}) = (f_j \circ \pi_i) \psi_{i1}(r_1) = (f_j \circ \pi_i) r_{i+1}.
\end{displaymath}
By Theorem 4.2 we conclude that 
\begin{displaymath}
\mathrm{dist} ((c_{i+1,1},c_{i+1,2},\ldots,c_{i+1,n}),
\mathrm{Lg}_n(C_{i+1})) \geq 1.
\end{displaymath}
As noted above, this implies that the simple limit $B_n$ has
stable rank strictly greater than $n$.

We now show that $\mathrm{sr}(B_n) \leq n+2$.
Given an exact sequence $B \rightarrow A \rightarrow C$ of
$C^*$-algebras, Corollary 4.12 of [9] states that
\begin{displaymath}
\mathrm{sr}(A) \leq \mathrm{max} \{ \mathrm{sr}(B), \mathrm{sr}(C)+1 \}.
\end{displaymath}
Applying this formula to the exact sequence $\mathrm{S}D_i 
\rightarrow A_i \rightarrow C_i$ we have
\begin{displaymath}
\mathrm{sr}(A_i) \leq \mathrm{max}\{
\mathrm{sr}(\mathrm{S}D_i),\mathrm{sr}(C_i)+1 \}.
\end{displaymath}
It is known that
\begin{displaymath}
\mathrm{sr}(p(\mathrm{C}(X) \otimes \mathcal{K})p) = \lceil \lceil \mathrm{dim}
  X / 2 \rceil / \mathrm{dim}p \rceil + 1 \nonumber
\end{displaymath}
whenever $X$ a compact Hausdorff space and $p$ is a projection in 
$\mathrm{C}(X) \otimes \mathcal{K}$ ([7]).  Thus, $\mathrm{sr} (C_i) = 
\mathrm{sr} (p_i (\mathrm{C}(X_i) \otimes \mathcal{K}) p_i) = n+1$ by 
inspection of the dimensions of the $p_i$ and $X_i$.  Since
$\mathrm{S}D_i$ is an ideal in $D_i \otimes \mathrm{C}([0,1])$, 
we have 
\begin{displaymath}
\mathrm{sr}(\mathrm{S}D_i) \leq \mathrm{sr}(D_i \otimes
\mathrm{C}([0,1])) \leq \mathrm{sr}(D_i) + 1
\end{displaymath} 
by Corollary 7.2 of [9].  
Theorem 6.1 of [9] states that 
\begin{displaymath}
\mathrm{sr}(\mathrm{M}_n (A))\leq \lceil (\mathrm{sr}(A)-1) / n \rceil
+ 1, 
\end{displaymath}
so that 
$\mathrm{sr}(D_i) = \mathrm{sr}(\mathrm{M}_{k_i \mathrm{dim}p_i}
\otimes C_i) \leq n+1$ for all $i$.
We conclude that $\mathrm{sr}(A_i) \leq n+2$, so that 
$\mathrm{sr}(B_n) \leq n+2$ by Theorem 5.1 of [9].
Combining this with the fact that $\mathrm{sr}(B_n) \geq n+1$ yields 
Lemma 4.1. \hfill $\Box$

\begin{lms}
The ordered group $\mathrm{K}_* B_n =\mathrm{K}_0 B_n \oplus \mathrm{K}_1
B_n$ is weakly unperforated.  Its order structure is the strict one 
coming from the first direct summand $(\mathrm{K}_0 B_n, {\mathrm{K}_0
  B_n}^+) = (\Z,{\Z}^+)$.
\end{lms}

\noindent
\textbf{Proof:}  
Since $(\mathrm{K}_0 B_n, {\mathrm{K}_0 B_n}^+)$ is weakly
unperforated it will be enough to show that every element in 
$\mathrm{K}_1 B_n$ is the $\mathrm{K}_1$-class of a
unitary element in $B_n$.  Since $\mathrm{K}_* B_n$ is the inductive
limit of the $\mathrm{K}_* A_i$, it will suffice to prove this assertion
for all $A_i$ with $i$ sufficiently large.  By the formulas and discussion in
the proof of Lemma 4.1, we know that $\mathrm{sr}(\mathrm{M}_{\mathrm{dim} p_i}(\mathrm{S}C_i))=2$ for 
all $i$ sufficiently large.  Assume that $i$ is so large for the remainder of the proof.   

From [1] and [9] we know that there is a bijection between elements of 
$\mathrm{K}_1 \mathrm{S}D_i$ and the $\mathrm{K}_1$-classes of unitaries in $\mathrm{M}_3 \otimes \widetilde{
\mathrm{M}_{\mathrm{dim} p_i}(\mathrm{S}C_i)}$.  Furthermore, any unitary in this latter
algebra is homotopic to a unitary in $\widetilde{\mathrm{M}_{3 \mathrm{dim} p_i}(\mathrm{S}C_i)}$.
Unitaries in $\widetilde{\mathrm{M}_{3 \mathrm{dim} p_i}(\mathrm{S}C_i)}$ give rise to 
unitaries in $\mathrm{S} D_i$, since $3 \leq k_i$ for all $i$.  Thus,
every element of $\mathrm{K}_1 \mathrm{S}D_i$ can be represented as the $\mathrm{K}_1$-class
of a unitary.  The map $\mathrm{K}_1 \iota$ induced by the inclusion
$\iota: \mathrm{S} D_i \rightarrow A_i$ is surjective (as
$\mathrm{K}_1 C_i = 0$) and
the desired conclusion for $A_i$ follows from functoriality. \hfill $\Box$

\begin{lms}
For $n > 2$, $B_n$ and $B_n \otimes  \mathcal{Z}$ are not isomorphic.
\end{lms}

\noindent
\textbf{Proof:}
We proceed by showing that $\mathrm{sr}(B_n \otimes \mathcal{Z}) \leq
2$, so that $\mathrm{sr}(B_n) \neq \mathrm{sr}(B_n \otimes \mathcal{Z})$. 

The algebra $\mathcal{Z}$ is an inductive limit of prime dimension drop algebras
 $\mathrm{I}[p_i,p_i q_i, q_i]$, $i = 1,2,\ldots$, where $p_i \rightarrow
\infty$ and $q_i \rightarrow \infty$ as $i \rightarrow \infty$ (cf. [6]).  For any $C^*$-algebra
$A$ the algebra $\mathrm{I}[p_i,p_i q_i, q_i] \otimes A$ is a full
algebra of operator fields, so by Theorem 1.1 of [8] we have
\begin{displaymath}
\mathrm{sr}(\mathrm{I}[p_i,p_i q_i, q_i] \otimes A) \leq \mathrm{sup}_{t \in [0,1]}
\{ \mathrm{sr}(A_t \otimes \mathrm{C}([0,1])) \},
\end{displaymath}
where $A_t$ is the fibre of $\mathrm{I}[p_i,p_i q_i, q_i] \otimes A$ at $t \in [0,1]$.
Since each such fibre is one of $\mathrm{M}_{p_i}(A)$, $\mathrm{M}_{q_i}(A)$, or 
$\mathrm{M}_{p_i q_i}(A)$ we may rewrite our estimate above as
\begin{eqnarray}
\mathrm{sr}(\mathrm{I}[p_i,p_i q_i, q_i] \otimes A) & \leq & \mathrm{max}
\{ \mathrm{sr}(\mathrm{M}_{p_i q_i}(A \otimes \mathrm{C}([0,1]))),
\nonumber \\
& & \ \ \mathrm{sr}(\mathrm{M}_{q_i}(A \otimes \mathrm{C}([0,1]))),
\mathrm{sr}(\mathrm{M}_{p_i}(A \otimes \mathrm{C}([0,1]))) \}. \nonumber 
\end{eqnarray}
By Corollary 7.2 of [9] we have $\mathrm{sr}
(A \otimes C[0,1]) \leq \mathrm{sr}(A) + 1$.  By Theorem 6.1 of [9] we have 
that $\mathrm{sr}(\mathrm{M}_n (A))
\leq \lceil (\mathrm{sr}(A)-1) / n \rceil + 1$.  Thus, there exists $i_0 \in
\mathrm{N}$ such that $\mathrm{sr}(\mathrm{M}_{p_i q_i}(A \otimes \mathrm{C}([0,1])))$,
$\mathrm{sr}(\mathrm{M}_{q_i}(A \otimes \mathrm{C}([0,1])))$ and
\mbox{$\mathrm{sr}(\mathrm{M}_{p_i}(A \otimes \mathrm{C}([0,1])))$} are all less than
or equal to two for $i \geq i_0$.  We conclude that 
\begin{displaymath}
\mathrm{sr}(\mathrm{I}[p_i,p_i q_i, q_i] \otimes A) \leq 2
\end{displaymath}
for all $i \geq i_0$.  Finally, $B_n \otimes \mathcal{Z}$ is an
inductive limit of algebras of the form $\mathrm{I}[p_i,p_i q_i, q_i] \otimes B_n$, 
all but finitely many of which have stable rank less than or equal to
two.  By Theorem 5.1 of [9], the limit $B_n \otimes  \mathcal{Z}$ must
have stable rank less than or equal to two, as claimed. \hfill $\Box$

Thus, we have established Theorem 1.1.  In closing, we note that given
two natural numbers $n$ and $m$ one may carry out the construction of
section 3 to produce algebras $B_n$ and $B_m$ which, if the parameters
$q_i$ are chosen to be the same for both constructions, will have
isomorphic Elliott invariants.  This shows that one can produce
simple, nuclear, infinite-dimensional, stably finite counterexamples
to the Elliott conjecture which lie entirely outside the class of
$\mathcal{Z}$ absorbing $C^*$-algebras.  The explicit calculation of
$\mathrm{Ell}(B_n)$ and $\mathrm{Ell}(B_m)$ is long and not
particularly illuminating.  We leave it to the reader.

\pagebreak

\noindent
Institute for Mathematical Sciences, K{\o}benhavns Universitet,
Universitetsparken 5, DK-2100 K{\o}benhavn {\O}, Danmark \newline
Electronic mail: atoms@math.ku.dk

\end{document}